\documentclass{article}

\long\def\ignore#1{\relax}

\usepackage{a4}
\usepackage{amssymb}
\usepackage{QED}
\usepackage{url}

\newtheorem{theorem}{Theorem}

\newtheorem{corollary}[theorem]{Corollary}

\newtheorem{definition}[theorem]{Definition}

\newtheorem{lemma}[theorem]{Lemma}
\newtheorem{notation}[theorem]{Notation}

\newtheorem{proposition}[theorem]{Proposition}
\newtheorem{remark}[theorem]{Remark}

\def\l{\lambda}

\def\f{\rightarrow}
\def\fst{\f_{st}}

\def\<{\langle}
\def\>{\rangle}

\def\R{\ifmmode{\rm I\mkern-3.1mu R\mkern1mu}\else{\rm
I\kern-.18em  R\hskip1pt\ }\fi\relax}

\def\Z{\ifmmode{ Z\mkern-8.0mu Z\mkern2mu}\else{
Z\kern-.32em Z\hskip1pt\ }\fi\relax}

\def\Q{\ifmmode{\rm Q\mkern-10mu
l\mkern4.5mu}\else{\rm Q\kern-.57em l\hskip3pt\
}\fi\relax}

\def\N{\ifmmode{\rm I\mkern-3.1mu
N\mkern0.5mu}\else{\rm I\kern-.16em N\hskip0.5pt\
}\fi\relax}

\def\C{\ifmmode{\rm C\mkern-8.8mu l\mkern4mu}\else{\rm
C\kern-.48em l\hskip2.6pt\ }\fi\relax}

\def\tr{\triangleright}

\def\ras{\tr^*}
\def\ih{{\em IH}}
\def\ihb{{\em IH }}

\newcommand{\bra}[1]{[#1]}
\def\sig{\sigma}

\def\notin{\not\in}
\title{A direct proof of   the  confluence \\ of  combinatory strong reduction}

\author{Ren\'e David  \\ Universit\'e de Savoie, Campus
Scientifique\\
73376 Le Bourget du Lac, France. \\
 Email : david
@univ-savoie.fr}

\begin{document}
\maketitle \pagestyle{plain}

\begin{abstract}
I give a proof of  the  confluence of  combinatory strong
reduction that does not use the one of $\lambda$-calculus. I also
give simple and direct proofs of a standardization theorem for
this reduction and the strong normalization of simply typed terms.

\end{abstract}
\section{Introduction}

 Combinatory Logic (see \cite{curry}, \cite{curry&all}) is a
first order language that simulates the $\lambda$-calculus without
using bounded variables. \ignore{It has the same properties as
 the
$\lambda$-calculus: confluence, strong normalization of typed
terms, etc.} But, at present, the known proofs of confluence are
all based on the confluence of the $\lambda$-calculus which  has
to be proved before and thus Combinatory Logic is not a
self-contained theory. The question of getting a direct proof of
this confluence was raised long ago in \cite{curry} and appears in
the TLCA list of open problems. I give here such a proof.

 The paper is organized as follows. Section \ref{s2} gives the main
 definitions of Combinatory Logic, states the theorem  and the idea of the proof. Section \ref{s3} gives the proof of
 the confluence of an auxiliary system. Section \ref{s4} gives the equivalence of the two systems and
 deduce the confluence of the original one. Section \ref{sst} gives a standardization theorem
 and section \ref{ssn} gives a direct proof of strong normalization for simply typed terms.
 Finally, I conclude in section \ref{s5} with some remarks.

\section{The idea of the proof of confluence}\label{s2}
\subsection{Combinatory Logic}

\begin{definition}

 The set $C$ of combinators is defined by the following grammar
(where $x$ denotes a variable)

$$C:= x  \; \mid \; K  \; \mid \;  S  \; \mid \;  I \; \mid  (C \; C)$$
\end{definition}

In the literature, the objects determined by this grammar are
usually
 called CL-terms and the word combinator is given for
{\em closed} CL-terms. However since, in section \ref{s3},   the
word term will be used for something slightly different, I prefer
to keep the word combinator here.

\begin{definition}\label{def_orig}
 For $u \in C$, the term $[x] u$ is defined, by induction on $u$,
  by the following rules

 \begin{enumerate}
 \item $\bra{x} u= K
u$ if $x \not\in u$
\item $\bra{x} x = I$
\item $\bra{x} (u \ x) =u$ if $x \not\in u$
\item $\bra{x} (u \; v) = (S \; \bra{x} u \; \bra{x} v)$ if none of
the previous rules apply.
  \end{enumerate}

\end{definition}

\begin{definition}\label{red_orig}
 The reduction on combinators is the closure by contexts of the
  following rules.
\begin{enumerate}
  \item $(K \; u \; v )\succ u$ \hspace{1cm} $(S \; u\; v \; w) \succ  (u \; w \; (v \;
    w))$ \hspace{1cm} $(I \ u) \succ  u$
  \item $[x]u   \succ [x] v$  \hspace{1cm}  if $u \succ  v$
\end{enumerate}

\end{definition}

I recall here usual notions about reductions.

\begin{definition}
Let $\rightarrow$ be a notion of reduction.
\begin{itemize}
\item As usual,  $\rightarrow^*$ denotes the reflexive and transitive
  closure of $\rightarrow$.
  \item The reduction  $\rightarrow$ is locally confluent if, for any term
  $u$, the following holds. If $u \rightarrow u_1$ and $u
\rightarrow u_2$, then $u_1 \rightarrow^* u_3$ and $u_2
\rightarrow^* u_3$ for some $u_3$.

\item The reduction $\rightarrow$ commutes with the reduction $\rightarrow_1$ if, for any term
  $u$, the following holds. If $u \rightarrow^* u_1$ and $u
  \rightarrow_1^*
u_2$ then $u_1 \rightarrow_1^* u_3$, $u_2 \rightarrow^* u_3$ for
some $u_3$
  \item The reduction  $\rightarrow$ is  confluent  if it commutes
  with itself.
\item A  term $u$ is strongly
normalizing (denoted as $ u
  \in SN$) if there is no infinite reduction of $u$.
\end{itemize}

\end{definition}

\begin{remark}
Rule (2) of   definition \ref{red_orig} is fundamental to have the
equivalence of
  combinatory logic (denote as $LC$)
   and $\l$-calculus (denoted as $\Lambda$) in the following
   sense.  Let $H$ be the translation between
  $\Lambda$ and $LC$ defined by

\begin{center}
$H(x)=x$ \hspace{1cm} $H((u_1 \ u_2))=(H(u_1) \ H(u_2))$
\hspace{1cm}
  $H(\l x. u)= [x] H(u)$
\end{center}

\noindent Without rule (2), the compatibility property between
$\Lambda$ and $LC$ (i.e. if $t$ reduces to
   $t'$, then $H(t)$ reduces to $H(t')$)
  would not be true. This is because the reduction in $LC$ will
  not allow a reduction below a $\l$. For example, let $t=\l x.
  (\l y. x \ x)$. Then $H(t)=[x](K \ x \ x)=(S \ K \ I)$ is normal
  whereas $t$ is not.

\end{remark}

 Note that without rule (2) of   definition \ref{red_orig} (this reduction is then called weak reduction), the confluence would be trivially proved
  by the method of parallel reductions.

\begin{remark}\label{rq2}
The confluence of the reduction $\succ$ depends on the good
interaction between  rule (2) of   definition \ref{red_orig} and
clause (3) of definition \ref{def_orig} (which corresponds,
intuitively, to the $\eta$-equality of the $\lambda$-calculus). In
fact, the confluence of $\succ$ would not be true if clause (3) of
definition \ref{def_orig} had been omitted. The reason is the
following. Let $u$ and $v$ be two combinators. Assume variable $x$
occurs in $u$ but not in $v$ and $u$ reduces to $u'$ for some $u'$
that does not contain $x$ (for example $u=K \ y \ x$). Then, by
applying rule (2) of   definition \ref{red_orig}, we have \\
(1)\hspace{1cm}  $[x](u \ v)=(S \ [x]u \ (K \ v)) \succ (S \ [x]u'
\ (K \ v))= (S \ ( K \ u') \ (K \ v))$

 and
\\(2) \hspace{1cm} $[x](u \ v) \succ [x](u'  \ v) =(K \ (u'  \ v))$

\begin{itemize}
  \item Without clause (3) of definition \ref{def_orig} the term $(S \ ( K \ u') \ (K \ v))$
  is not of the form $[x]w$, so that the two terms $(K \ (u'  \
  v))$ and $(S \ ( K \ u') \ (K \ v))$ are not reducible to a common term.
  \item With clause (3) of definition \ref{def_orig} the term $(S \ ( K \ u') \ (K \ v))$
  is  of the form $[x]w$.

  \[ (S \ ( K \ u') \ (K \ v))= [x](S \ ( K \ u') \ (K \ v) \ x)  \]

  from which

  \[ [x](S \ ( K \ u') \ (K \ v) \ x) \succ [x](K \ u' \ x \ (K \ v \ x )) \succ^* [x](u' \ v) = (K \ (u'  \
  v)) \]

 \noindent  Thus the two terms $(K \ (u'  \
  v))$ and $(S \ ( K \ u') \ (K \ v))$ are  reducible to a common term.

\end{itemize}

\end{remark}

The main result of this paper is the following theorem.

\begin{theorem}\label{thm}
The reduction $\succ$ on combinators is confluent.
\end{theorem}

\subsection{The idea of the proof}\label{s2.2}

 I want to prove the confluence  by using the same method
as in \cite{david} i.e. by proving first a theorem on finiteness
of developments. Then, by this theorem, Newman's Lemma and the
local confluence of the developments  we get the confluence of
developments. Then it remains to show that the reduction itself
 is the transitive closure of the developments.

 But the given
system is quite hard to study because it is difficult to mark the
redexes and thus to give a precise definition for a theorem on
finiteness of developments. This is also because the form of a
term does not determine easily its redexes. The main technical
reason is the following. We should think that any reduct of $[x]u$
would have the form $[x]u'$ for some reduct $u'$ of $u$. But this
property, which is trivial in the $\lambda$-calculus, is not true
here. Here is an example. Let $u,v$ be combinators, $x$ be a
variable that occurs both in $u$ and $v$ and let $t=[x](u \ v)=(S
\ [x]u \ [x]v)$. Then, it is easy to check that $t=([y][x](u \ (y
\ x)) \ [x]v)$. Now if $u=(S \ u_1 \ u_2)$ then $t$ reduces to
$t'=([y][x](u_1 \ (y \ x) \ (u_2 \ (y \ x))) \ [x]v)$ and it is
easy to check that $t'$ cannot be written as $[x]w$ for some
reduct $w$ of $(u \ v)$. Note that, in the $\lambda$-calculus, the
corresponding  equality i.e. $\lambda x. (u \ v)=(\lambda yx. (u \
(y \ x))\ \lambda x.v)$ needs $\beta$-reductions and not only
$\eta$-reductions whereas in Combinatory Logic it only comes from
the $\eta$-rule.

Thus I will first prove the confluence of an auxiliary system.
This system will be shown to be equivalent to the other one in the
sense  that the symmetric and transitive closure of both systems
are the same. Then I will deduce the confluence of the first
system from the one of the second.

The auxiliary system treats separately the reductions that,
intuitively,   corresponds in the $\lambda$-calculus to $\beta$
and $\eta$. To prove the confluence of this system, I prove the
confluence of $\beta$. This is done, as mentioned above, by
proving a theorem on finiteness of developments. Note that the
fact that the reduction is the transitive closure of developments
 (which is trivial in the $\lambda$-calculus) is not so easy here.  I
deduce the confluence of the whole system (intuitively $\beta$ and
$\eta$) by another commutation lemma.

\begin{lemma}[Newman's Lemma]\label{new}
Let $\rightarrow$ be a notion of reduction that is locally
confluent and strongly normalizing. Then $\rightarrow$ is
confluent.
\end{lemma}

\section {An auxiliary system}\label{s3}

To define this new system, I first remove the $\eta$-equality in
the definition of the abstraction.

\begin{definition}\label{def2}
\begin{enumerate}
 \item $\lambda x. u= (K
\ u)$ if $x \not\in u$
\item $\lambda x. x = I$
\item $\lambda x. (u \; v) = (S \; \lambda x. u \; \lambda x. v)$ if none of
the previous rules apply.
  \end{enumerate}
\end{definition}

\noindent and I  add new reduction rules. In Definition
\ref{def_red} below rule (2) is necessary to have confluence. Rule
(3) corresponds to the $\eta$-reduction and is necessary to have
the equivalence with the other system.

\begin{definition}\label{def_red}
\begin{enumerate}
  \item $(K \; u \; v )\rightarrow u$ \hspace{.5cm} $(S \; u\; v \; w) \rightarrow (u \; w \; (v \;
    w))$ \hspace{.5cm} $(I \ u) \rightarrow u$
    \item $(S \ (K \ u) \ (K \ v)) \rightarrow (K \ (u \ v))$
    \item $(S \ (K \ u) \ I) \rightarrow u$
  \item $\lambda x.u   \rightarrow \lambda x. v$  \hspace{1cm}  if $u \rightarrow v$
\end{enumerate}

\end{definition}

It is important to note that the two reductions $\succ$ and
$\rightarrow$ are not the same i.e. there are combinators such
that $u \rightarrow^* v$ for some $v$ but $u$ does not reduce to
$v$ by  $\succ$ and, similarly,  there are combinators such that
$u \succ^* v$ for some $v$ but $u$ does not reduce to $v$ by
$\rightarrow$. Here are examples. Let $u=[y][x](S \ x \ x \ (y \
x))$. Then $u \succ [y][x](x \ (y \ x)\  (x \ (y \ x)))$ and it is
easy to check that $u$ is normal for $\rightarrow$. Let
$u_1=\lambda x. (S \ y \ x \ x) \rightarrow \lambda x. (y \ x \ (x
\ x))=v$ and it is not too difficult to check that $u$ does not
reduce to $v$ by  $\succ$.

Although the two reductions $\succ$ and $\rightarrow$ are not the
same, we now show that they give the same equations on
combinators. I denote by $\equiv$ the equivalence relation induced
by $\succ$ i.e. $u \equiv v$ iff  there is a sequence $u_0, ...,
u_n$ of combinators such that $u_0=u$, $u_n=v$ and, for each $i$,
either $u_i \succ u_{i+1}$ or $u_{i+1} \succ u_i$. The equivalence
induced by $\rightarrow$ will be denoted by $\approx$.

\begin{lemma}
\begin{enumerate}
  \item For each $u,v$, $(S \ (K \ u) \ I) \succ^* u$ and $(S \ (K \
  u) \ (K \ v)) \succ^* (K \ (u \ v))$
  \item For each $u$, $\lambda x. u \rightarrow^*
[x]u$ and $\lambda x. u \succ^* [x]u$.

\end{enumerate}
\end{lemma}
\begin{proof}
\begin{enumerate}
  \item Let $x$ be a fresh variable. Then, $(S \ (K \ u) \ I)=[x](S \ (K \
u) \ I \ x) \succ [x] (K \ u \ x \ (I \ x)) \succ^* [x](u \ x ) =
u$ and $(S \ (K \
  u) \ (K \ v)) =[x] (S \ (K \
  u) \ (K \ v) \ x) \succ [x](K \ u \ x \ (K \ v \ x)) \succ^*
  [x](u \ v) = (K \ (u \ v)))$.
  \item This follows immediately from the first point.
\end{enumerate}
\end{proof}

\begin{theorem}\label{equiv}
Let $u,v$ be combinators. Then $u \equiv v$ iff $u \approx v$.
\end{theorem}
\begin{proof}
It is enough to show that if $u \succ v$ then $u \approx v$ and if
$u \rightarrow v$ then $u \equiv v$. Each  point is proved by
induction on the level of the reduction. The result is trivial for
the level 0. Assume then that the level is at least 1. For the
first direction, I have to show that, if $u \succ v$ then $[x]u
\approx [x]v$. By the \ih\ we know that $u \approx v$ and it is
thus enough to show that, if $u \rightarrow v$, then $[x]u \approx
[x]v$. By the previous lemma, we have $\lambda x. u \rightarrow
[x]u$ and, since $\lambda x. u \rightarrow \lambda x. v
\rightarrow [x]v$, we are done. For the other direction, we have
to prove that, if $u \succ v$, then $\lambda x. u \equiv \lambda
x. v$. This is because $\lambda x. u \succ [x] u \succ [x]v$ and
$\lambda x. v \succ [x]v$.
\end{proof}

\begin{theorem}\label{thm_aux}
The reduction $\rightarrow$ on combinators is confluent.
\end{theorem}

 As mentioned before, to prove this theorem I first prove the confluence of the system where
 the $\eta$-reduction (i.e. rule (3) of definition \ref{def_red}) has been removed. The theorem on finiteness of
developments of this system can be formalized as theorem
\ref{main} below. I need some new definitions.

\subsection{Some definitions}

\begin{definition}\label{lambda}
Let $V$ be an infinite set of variables.

\begin{itemize}
  \item Let  $A=V
\cup \{S_i \ / \ i=0,1,2,3\} \cup \{K_i \ / \ i=0,1 \}\cup \{I_i \
/ \ i=0,1 \} $. The elements of $A$ will be called   atoms.
  \item The set of  terms is defined by the following grammar

  $$T:= A \mid \   \ (T \; T)$$
  \item The size of a term (denoted as $size(t)$) is defined by the
  following rules: for $\alpha \in A$, $size(\alpha)=1$ and
  $size((u \ v))=size(u)+size(v)+1$.

\end{itemize}

\end{definition}

The meaning of the indices on $S, K, I$ is the following. First, I
want to mark the redexes that are allowed to be reduced. I do this
by simply indexing the letters $S, K, I$. The index 0 means that
the symbol is not marked (i.e. we are not allowed to reduce the
corresponding redex), the index 1 means that the redex
 is allowed.

I also want to indicate whether or not  a combinator $S, K, I$ is
the first symbol of  a term of the form $\lambda x.u$ for which I
want to reduce in $u$. Actually, for $K,I$ there is nothing to do
 because a variable has no redex and, since $\lambda x.u = (K \ u)$ when
 $x$ does not occur in $u$, the redexes in $u$ are,
 in fact, already visible at the top level.
But for $S$ this will be useful and I need thus 4 indices.

\begin{itemize}
  \item $S_0$ is an $S$ that is neither
marked nor introduced by a $\lambda$,
  \item $S_1$ is an $S$ that is
marked but not introduced by a $\lambda$,
  \item $S_2$ is an $S$ that is
not marked but introduced by a $\lambda$
\item $S_3$ is an
$S$ that is
 marked and introduced by a $\lambda$.
\end{itemize}

\begin{definition}\label{index-lambda}
Let $u$ be a term and $x$ be a variable. I define, for $i=0,1$ the
set of terms (denoted as $\lambda_i x. t$) by the following rules.
\begin{enumerate}
  \item if $t=x$, $\lambda_i x.  t=\{I_i\}$
  \item if $t \neq x$ is an atom,  $\lambda_i x. t=\{(K_i \ t)\}$
  \item if $t=(u \ v)$ and $x\not\in t$, $\lambda_i x. t=\{(K_i \ t)\} \cup \{(S_{i+2} \ u' \ v')\ | \ u' \in \lambda_i x.
u, v' \in \lambda_i x. v\}$
\item if $t=(u \ v)$ and $x \in t$,
$\lambda_i x. t=\{(S_{i+2} \ u' \ v')\ | \ u' \in \lambda_i x. u,
v' \in \lambda_i x. v\}$.
\end{enumerate}
\end{definition}

The reason  of this unusual definition and, in particular, the
fact that $\lambda_i x. t$ represents a set of terms instead of a
single term,  is the following.  It will be useful to ensure that
the set of terms of the form $\lambda x.u$ is closed by reduction.
But this is not true if the abstraction is defined by the rules of
definition \ref{def2}.

Here is an example. Let $u$ and $v$ be two combinators. Assume
variable $x$ occurs in $u$ but not in $v$ and $u$ reduces to $u'$
for some $u'$ that does not contain $x$. As shown in points (1)
and (2) of  Remark \ref{rq2},  $\lambda x.(u \ v)$ reduces to $(S
\ (K \ u') \ (K \ v))$ and $(K \ (u' \ v))$. Allowing, in such a
case, both $(K \ (u' \ v))$ and $(S \ (K \ u') \ (K \ v))$ to be
in $\lambda x.(u' \ v)$ will repair this problem.

The given definition  is then an indexed version of this idea. The
index 1 (resp. 0) will mean that the $S, K, I$ introduced by the
definition are marked (resp. are not marked) and thus allow a
redex to be reduced. Note that the $i+2$ indexing $S$ means
(depending whether $i=1$ or $i=0$) that $S$ comes from a $\lambda$
and is (or is not) marked.

\begin{definition}\label{reduction}
The reduction (denoted as $t \tr t'$) on  terms is the closure by
contexts of the following rules
\begin{enumerate}
  \item
  \begin{enumerate}
  \item For $i=1,3$ $(S_i \; u\; v \; w) \tr (u \ w \ (v \ w))$
  \item $(K_1 \ u \ v) \tr u$  and   $(I_1 \ u) \tr u$
  \item For $i=0,1$   $(S_{i+2} \ (K_i \ u) \ (K_i \ v))
  \tr (K_i \ (u \ v))$
  \ignore{\item  $(S_3 \ (K_1 \
  u) \ I_1) \tr u$
  \item   $(S_2 \ (K_0 \
  t) \ I_0) \tr t$  if  $t$ is non equal to either $I_1$ or
  $(K_1 \ u)$ or  $(S_1 \ u \ v)$.
  \item $(S_2 \ (K_0 \
  I_1) \ I_0) \tr I_0$, $(S_2 \ (K_0 \
  (K_1 \ u) \ I_0) \tr (K_0 \ u)$ and $(S_2 \ (K_0 \
  (S_1 \ u \ v) \ I_0) \tr (S_0 \ u \ v)$}

  \item For $i=0,1$, if $u \tr v$, $t \in \lambda_i x.  u$ and $ t' \in  \lambda_i x.
  v$, then $t \tr t'$
\end{enumerate}
  \item The level of a reduction
  (denoted as $lvl(t \tr t')$) is defined as follows.

\begin{itemize}
  \item If $t \tr t'$ by using   rule (a),(b) or (c), the level is 0.
  \item $t \tr t'$ by using   rule (d), the level is $lvl(u \tr v)+1$.

\end{itemize}

\end{enumerate}

\end{definition}

\noindent {\bf Remarks and examples}

These rules correspond to the indexed version of the rules (1, 2,
4) of definition~\ref{def_red} combined with the fact that
$\lambda x. u$ now is a set of terms.

 For example, if $x$ does not
occur in $(u \ v)$ and $u \tr u'$, since $(K_i \ (u \ v)) \in
\lambda_i x. (u \ v)$ and $(S_{i+2} \ (K_i \ u') \ (K_i \ v)) \in
\lambda_i x. (u' \ v)$ we have $(K_i \ (u \ v)) \tr (S_{i+2} \
(K_i \ u') \ (K_i \ v))$. Note that $(K \ (u \ v))$ does not
reduce to  $(S \ (K \ u') \ (K \ v))$ by the rules of
Definition~\ref{def_red}.


\subsection{Fair terms}

We will show the confluence of $\tr$ not of the entire set of
terms but on some subset (the set of fair terms) that we now
define. This is because we need a set that is closed by reduction
(see Lemma \ref{sigma}).

\begin{notation}
\begin{itemize}
  \item Let $E$ be a set of terms and $\overrightarrow{u}$ be a sequence
of terms (resp. $f$ be  function into terms). I will write
$\overrightarrow{u} \in E$ (resp. $f \in E$) to express the fact
that each term of the sequence $\overrightarrow{u}$ (resp. in the
image of $f$) is in $E$.
  \item Let $\overrightarrow{u}$ be  a finite (possibly empty) sequence
of terms and $v$ be a term.
 I denote by  $(v \ \overrightarrow{u})$  the term $(v \ u_1 \  ... \ u_n)$
   where $\overrightarrow{u}= u_1, ..., u_n$.
\end{itemize}
\end{notation}

\begin{definition}
\begin{itemize}
  \item An address is a finite list of elements of the set
  $\{l,r\}$.
  \item The empty list will be denoted by $\varepsilon$ and
  $[a::l]$ (resp. $[l::a]$) will denote the list obtained from $a$ by
  adding $l$  at the end (resp. at the beginning) of $a$ and
  similarly for $r$.
  \item If $a,a'$ are addresses, I will denote by $a < a'$ the
  fact that $a$ is an initial segment of $a'$.
  \item  Let $u$ be a term. I will denote by $u_a$ the sub-term of
  $u$ at the address $a$. More precisely, $u_a$ is defined by the following rules: $u_\varepsilon=u$, $(u \
  v)_{[l::a]}=u_a$ and $(u \
  v)_{[r::a]}=v_a$.
\end{itemize}
\end{definition}

\begin{definition}\label{def-phi}
\begin{itemize}

  \item Let $u$ be a term and $f$ be a function from a set $E$ of addresses in $u$ into terms. I say that
  $f$ is adequate for $u$ (I will also say $(u,f)$ is  adequate) if there are no addresses $a,a'$ in $E$ such that $a <
  a'$.

  \item Let $(u,f)$ be  adequate and
   $x$ be a variable. Then $\phi_x(u,f)$ is
  a term obtained by replacing in $u$, for each $a\in dom(f)$,  the term at address
  $a$  by $(w_a \ f(a))$ for some $w_a \in \lambda_1 x.u_a$.
  \item Let $u$ be a term, $x_1, ..., x_n$ (resp. $f_1, ..., f_n$) be a
sequence (possibly empty) of variables (resp.  of functions). The
term $\phi_{x_1}(\phi_{x_2}(...(\phi_{x_n}(u, f_n),f_{n-1})
...)f_1)$ will be denoted by $\phi(u, \overrightarrow{x},
\overrightarrow{f})$ or simply $\phi(u)$ if we do not need to
mention explicitly  $\overrightarrow{x},\overrightarrow{f}$ or if
they are clear from the context.

\end{itemize}
\end{definition}

\noindent {\bf Comments and examples}

A typical term of the form $\phi_x(u,f)$ is obtained as follows.
Let $t= (\lambda_1 x.u \ v) $. First  reduce the head redex of $t$
(this intuitively means: do the $\beta$-reduction and introduce a
kind of explicit substitution $[x:=v]$) and then propagate (not
necessarily completely)  this substitution inside $u$ (this
intuitively means do some $S,K,I$ reductions at the top level),
possibly doing some (different) reductions in the (different)
occurrences of $v$.  The term obtained in this way is a typical
term of the form $\phi_x(u,f)$. Here is an example.

 Let $u=(y \ x \ x)$, $v, v'$ be combinators and let $f$ be such that
 $f([l])=v$ and $f([r])=v'$. Then $\phi_x(u,f)= (S_3 \ (K_1 \ y) \ I_1 \ v \ (I_1 \
 v'))$. Remark that, if  $v \tr
 v'$, we have $(\lambda_1 x. u \ v) \tr \phi_x(u,f)$.

 Note that, even if we only need $\phi_x(u,f)$ in case the terms in the image
 of $f$ are reducts of a single term, we do not ask this property in the
 definition.

 Finally note that, in the same way that $x$ does not occur in $\lambda_ix. u$, it does not occur in
 $\phi_x(u,f)$. This implies that, as usual, when we substitute a variable
 $y$ by some term $v$ in a term of the form  $\lambda_ix. u$ or
 $\phi_x(u,f)$  we may assume (by possibly renaming $x$ with a fresh name) that $x$ does not occur in
 $v$,
avoiding then its capture.

\begin{definition}\label{fair}
The set $F$ of fair terms is defined by the following grammar.
\begin{enumerate}
  \item $x, S_0, K_0, I_0$ are fair
  \item If $u,v$ are fair  then so is $(u \ v)$.
  \item If $u$ is fair and $t \in \lambda_0  x. u$ then so is $t$.
  \item If $v_1, v_2, v_3$ are fair, then so are $ (S_1 \ v_1 \
v_2 \  v_3)$, $(K_1 \ v_1 \  v_2 )$ and $ (I_1 \ v_1)$
\item If $x$ is a variable, $u,f \in F$ and $(u,f)$ is  adequate,  then
$\phi_x(u,f)$ is fair.
\end{enumerate}
\end{definition}

\ignore{Note that, in rule (3),  I have written $t \in \lambda_0
x. u$ instead of $t = \lambda_0  x. u$ to remind that $\lambda_0
x. u$ represents a family of terms. In the rest of the paper, when
the expression $\lambda_i  x. u$ appears somewhere, this means
that it can be interpreted as any of the terms of the family. I
will write only $t \in \lambda_0  x. u$ when I want to focus on
the fact that it is a family.}

 Fair terms are thus combinators where we have marked the redexes that are allowed to be reduced.
 The terms of the form $\phi_x(u,f)$
are introduced for the following reason. If $t=(w \ v)$ for some
$w \in \lambda_1  x.u $, I may want to reduce both a redex in $u$
and $t$ as a redex. Thus the set of fair terms must be closed by
the following rule: (6) If $u, v$ are fair then so is $t=(w \ v)$
for $w \in \lambda_1 x.u$. But, if I had defined fair terms by
rules 1, 2, 3, 4 and 6, then $F$ will not be closed by reduction
because, if $w \in \lambda_1  x. u$, the reduct of $t=(w \ v)$
will not necessarily be fair. The reason is the following. Let $u=
(u_1 \ u_2)$ be such that $u$ is fair but $u_1$ is not (for
example $u_1=(K_1 \ y), u_2=y$). Then $v=(\l_1  x.  u \ z)$ is
fair. But $v \tr v'= (\l_1 x. u_1 \; z) \ (\l_1 x. u_2 \; z)$ and
$v'$ may not be fair since $u_1$ is not.

\begin{definition}
Let $u$ be fair. I denote by $nb(u)$ the number of rules that have
been used to prove that $u$ is fair.
\end{definition}

\subsection{Some properties of fair terms}

\begin{lemma}\label{subs2}
The set of fair terms is closed by substitutions.
\end{lemma}
\begin{proof}
 By an immediate induction on $nb(u)$. Use
 the fact that, if  $t \in \lambda_i  x. u$, then $\sig(t) \in  \lambda_i
 x.\sig(u)$.
\end{proof}

\begin{lemma}\label{carac}
Let $t=(\alpha \ \overrightarrow{u})$ be fair
 where $\alpha$ is an atom.

\begin{enumerate}
  \item If $\alpha$ is $S_2$, then
$lg(\overrightarrow{u}) \geq 2$. If $\alpha$ is $S_1$ or $S_3$,
then $lg(\overrightarrow{u}) \geq 3$.
  \item If $\alpha$ is $K_1$, then
$lg(\overrightarrow{u}) \geq 2$. If $\alpha$ is $I_1$, then
$lg(\overrightarrow{u}) \geq 1$.
\end{enumerate}

\end{lemma}
\begin{proof}
By induction on $nb(t)$.
 I only look at the cases with $S$. The other ones are similar.

\begin{itemize}
  \item If
  the last rule that has been used to prove $t \in F$ is  (2) of definition
  \ref{fair},
  the result follows immediately from the \ih. If it is rule (4) the result is
  trivial.

  \item If it is rule (3). If $\alpha=S_2$, the result is also trivial.
  The other cases are impossible.
  \item If it is rule (5) and $(\alpha \ \overrightarrow{u})=\phi_y(v,f)$.
   Let  $a$  be the leftmost address in $dom(f)$.
   For $\alpha=S_1$ (resp. $\alpha=S_2$)
  we may not have
   $a=[l,l, ...,l]$  since this will imply that $t$
   begins with $S_3$. Thus $v=(S_1 \
   \overrightarrow{w})$ (resp. $v=(S_2 \
   \overrightarrow{w})$) and the result follows from the \ih.
   For $\alpha=S_3$, if the leftmost address is not of the form  $[l,l,
   ...,l]$ the result is as before. Otherwise, this implies that
   $t = (w_a \ f(a) \
\overrightarrow{s})$ for some $w_a \in \lambda_1 y. v_a $ and some
$\overrightarrow{s}$ and the result is trivial.
  \end{itemize}
\end{proof}

\begin{lemma}\label{cara}
Let $u,u'$ be terms, $t\in \lambda_i y. u$ and  $t' \in \lambda_j
x. u'$. Assume $t$ is a sub-term of $t'$. Then, either $t$ is a
sub-term of $u'$ or $i=j$, $x=y$ and $u$ is a sub-term  of $u'$.
\end{lemma}
\begin{proof}
By induction on $u'$. \ignore{The only non trivial case is $u=(u_1
\ u_2)$. But it is easy to check that either $t= \lambda  x  \ u$
or $t$ is a sub-term of $\lambda  x  \  u_i$ for some $i$ and the
\ihb gives the result.}
\end{proof}

\begin{lemma}\label{carac1}
\begin{itemize}
  \item Let $t=(\alpha \
\overrightarrow{u}) \in F$ where $\alpha \in V \cup \{S_i, K_i,
I_i \ / i=0,1\}$. Then, $\overrightarrow{u} \in F$.
  \item If $t=(S_2 \ \overrightarrow{u}) \in F$, then $t=\phi((r \
\overrightarrow{w}))$ for some $r \in \lambda_0 y.v$ and some
$v,\overrightarrow{w} \in F$.
\end{itemize}
\end{lemma}

\begin{proof}
By  induction on $nb(t)$, essentially as in lemma \ref{carac}.
\end{proof}

\subsection{Some properties of reduction}

\begin{lemma}\label{decomp}
Let $ u_1,  u_2$ be fair and assume $t=(u_1 \  u_2) \tr t'$. Then
$t'=(u'_1 \  u_2)$ or  $t'=(u_1 \  u'_2)$ where $u_i \tr u'_i$.
\end{lemma}
\begin{proof}
It is enough to show that there is no possible
 interaction between $u_1$ and $u_2$. Such an interaction could occur in the following cases.
\\ - $lvl(t \tr t')=0$.  This is impossible  because, by Lemma \ref{carac}, all the arguments of the indexed
$S,K$ or $I$ of such a redex must be in $u_1$.
\\ - $lvl(t \tr t')>0$ and, for example,  $t\in  \lambda_0 x. v$ and $t'\in  \lambda_0 x. v'$ for some $v \tr v'$. This
could occur if
 $u_1=(S_2 \ w_1)$ for some $w_1 \in \lambda_0 x. t_1$, $u_2 \in  \lambda_0 x. t_2$ and $v= (t_1 \ t_2)$.
 But this is again impossible by Lemma \ref{carac}.
\end{proof}

\begin{lemma}\label{decomp1}
Let $u_1, u_2, u_3$ be terms.
\begin{itemize}
 \item Assume $t=(I_1 \ u_1)\tr t'$. Then either $t'=u_1$ or $t'= (I_1 \ u'_1)$ for $u_1 \tr u'_1$.
  \item Assume $t=(K_1 \ u_1 \ u_2) \tr t'$. Then either $t'=u_1$ or $t'=(K_1 \ u'_1 \ u_2)$ or
  $t'=(K_1 \ u_1 \ u'_2)$ for $u_i \tr u'_i$.
 \item Assume $t=(S_1 \ u_1 \ u_2 \ u_3) \tr t'$. Then either $t'=(u_1 \ u_3 \ (u_2 \ u_3))$ or
$t'=(S_1 \ u'_1 \ u'_2 \ u'_3)$ where $u_i \tr u'_i$ for a unique
$i$ and $u'_j=u_j$ for $j \neq i$. \ignore{ or $t=(S_1 \ (K \ v_1
\ I \ u_3)$ and $t'=(v_1 \ u_3)$ or $t=(S_1 \ (K \ v_1) \ (K \
v_2) \ u_3)$ and $t'=(K (v_1 \ v_2) \ u_3)$}
\end{itemize}
\end{lemma}

\begin{proof}
 It is enough to show that the mentioned reductions are the
 only possibilities. I only look at the last case since the other ones are similar.

 If  $ lvl(t \tr
 t')=0$, the result is trivial. Otherwise, this means that there
 is a sub-term of $t \in \lambda_i x. v$ which reduces
 to a term in $\lambda_i x.v'$  for $v \tr v'$. But, this sub-term
 has to be a sub-term of some $u_j$ because, otherwise (by Definition \ref{index-lambda})
 we will have $S_2$ or $S_3$ instead of
 $S_1$, and the result follows immediately.
\end{proof}

\begin{lemma}\label{lam}
Assume $t\in \lambda_0  x. u$ and $t \tr t'$. Then either $t'\in
\lambda_0 x. u$ and $size(t')< size(t)$ or $t'\in \lambda_0 x. u'$
 for some $u'$ such that $u \tr u'$.
\end{lemma}
\begin{proof}
If $lvl(t \tr t')=0$, the reduction cannot use (the closure by
context of) rule (a) in Definition \ref{reduction}. This is
because, since $t\in \lambda_0  x. u$, the index of $S$ in the
reduced redex cannot be 1 or 3 and thus the result is clear.
Otherwise, this follows easily from Lemma \ref{cara}.
\end{proof}

\begin{lemma}\label{carac5}
Assume $\phi (u, \overrightarrow{y}, \overrightarrow{f})\in
\lambda_0 x.v$. Then $u\in \lambda_0 x.w$ for some $w$ such that
$\phi(w, \overrightarrow{y}, \overrightarrow{f})=v$.
\end{lemma}
\begin{proof}
By an immediate induction on the length of the sequence
$\overrightarrow{y}$ it is enough to prove the result for $\phi_y
(u, f)$. This is proved by induction on $v$. I only consider the
case $v=(v_1 \ v_2)$ and $\phi_y (u, f)=(S_2 \ r_1 \ r_2)$ where
$r_j \in \lambda_0 x.v_j$ (the other cases are similar). The
leftmost address in $dom(f)$ cannot be $[l,l,...,l]$ because,
otherwise, $\phi_y (u, f)$ will begin with $S_3$. Thus $u$ is an
application  and $\phi_y (u, f)=(S_2 \ \phi_y (u_1, f_1) \ \phi_y
(u_2, f_2))$ where $u =(u_1 \ u_2)$. Thus $\phi_y (u_i, f_i) \in
\lambda_0 x.v_i$ and we conclude by the \ih.
\end{proof}

\begin{lemma}\label {phi}
Let $u,f \in F$ be such that $(u,f)$ is adequate. Then a redex in
$t=\phi_x(u,f)$ is either  in $u$ or in some $f(a)$ or is $(w_a \
f(a))$ for some $a$ and some $w_a \in \lambda_1 x. u_a$. Thus,  if
$t \tr t'$, one of the following cases holds.
\begin{itemize}
  \item $t'=\phi_x(u',f')$ for some $u',f'$ such that $u \tr u'$
  \item $t'=\phi_x(u,f')$ where $f \tr f'$
  \item $t'$ is obtained from $t$ by reducing the redex $(w_a \ f(a))$ for some
  $a\in dom(f)$ and some $w_a \in \lambda_1 x. u_a$. Then, $t'=\phi_x(u', f')$ and

\begin{itemize}
  \item If $u_a=x$, then $u'$ is $u$ where the occurrence of $x$ at the
   address $a$ has been replaced by $f(a)$ and  $dom(f')=dom(f) - \{a\}$.

  \item If $x \notin u_a$, then $u'=u$  and  $dom(f')=dom(f) - \{a\}$.

\item If $u_a=(v_1 \ v_2)$ then $u'=u$,  $dom(f')=dom(f)-\{a\} \cup
\{[a::l],[a::r]\}$, $f'([a::l])=f'([a::r])=f(a)$ and, for $b \neq
a, f'(b)=f(b)$.
\end{itemize}
\end{itemize}
\end{lemma}
\begin{proof}
By induction on $nb(u)$. The only thing to be shown is that the
mentioned cases are the only possible ones.  For $lvl(t \tr
t')=0$, this follows immediately from the fact that
terms of the form $(w_a \ f(a))$ for some $w_a \in \lambda_1 x.
u_a$  cannot introduce an interaction since they are redexes. For
$lvl(t \tr t')
> 0$, assume $r \in \lambda_i x. w$ is  a sub-term of $ \phi_y(u,f)$  and the reduction takes places in $w$.
Then, by Lemma \ref{carac5}, either the reduction is actually in
$f$ or  $w=\phi_y(v',f')$ for some adequate $(v',f')$ and the
result follows from the \ih.
\end{proof}

\begin{lemma}\label{sigma}
\begin{itemize}
  \item The set of fair terms is closed by reduction.
  \item Let $u$ be fair and $\sig$ be a fair substitution. Assume  $t=\sig(u) \tr t'$, then
either $t'=\sig(u')$ for some $u \tr u'$ or $t'=\sig'(u)$ for some
$\sig \tr \sig'$.

\end{itemize}
\end{lemma}
\begin{proof}
By induction on $nb(u)$, using Lemmas \ref{decomp}, \ref{decomp1},
\ref{lam} and \ref{phi}.
\end{proof}

\subsection{Confluence of $\tr$ on  fair terms}

\begin{lemma}\label{subs1}
Let $u$ be fair and $\sig$ be a fair substitution. If  $u, \sig
\in SN $, then so is $\sig(u)$.
\end{lemma}

\begin{proof}
This follows immediately from Lemma \ref{sigma}.
\end{proof}

\begin{theorem}\label{main}
Any fair term $t$ is in $SN$.
\end{theorem}

\begin{proof}
By induction on $nb(t)$.

\begin{itemize}
  \item If $t=x, S_0, K_0, I_0$, the result is trivial.
  \item If $t=(t_1 \ t_2)$, then, by the \ih,
$t_1,t_2 \in SN$ and, since $t=\sig((x \ y))$ where $\sig(x)=t_1$
and $\sig(y)=t_2$, the result follows from Lemma \ref{subs1}.
  \item If $t=(S_1 \ t_1 \ t_2 \ t_3)$,
$t= (K_1 \ t_1 \ t_2)$ or $t= (I_1 \ t_1)$ the proof is similar,
e.g. $(S_1 \ t_1 \ t_2 \ t_3)=\sig((S_1 \ x_1  \ x_2 \ x_3)$ where
$\sig(x_i)=t_i$.
\item If $t\in \lambda_0  x.  v$, the result follows from Lemma \ref{lam} and the \ih.
\item  Finally,
assume $t=\phi_x(u, f)$. Let $t'$ be the term obtained from $u$ by
replacing, for each $a \in dom(f)$, $u_a$ by $u_a[x:=f(a)]$. It
follows from Lemma \ref{subs1} that $t' \in SN$. But, by Lemma
\ref{phi}, and infinite reduction of $t$ would give  an infinite
reduction of $t'$ since it is not possible to have infinitely many
successive reductions of $t$ of the form of the last case of Lemma
\ref{phi}. Thus $t$ is in $SN$.
\end{itemize}
\end{proof}

\ignore{
\begin{lemma}\label{exemple}
Let $u,v$ be terms. Then $(S \ (K \ u) \ I) \tr u$, $(S \ (K \ u) \ (K \ v)) \tr (K \ (u \ v))$,
 $(S_1 \ (K_1 \ u) \ I_1) \tr u$ and $(S_1 \ (K_1 \ u) \ (K_1 \ v)) \tr (K_1 \ (u \ v))$.
\end{lemma}

\begin{proof}
and \\
 The proof is similar for the two other reductions.
\end{proof}
}

\begin{lemma}
Let $u,v$ be terms. Then, for $w \in \lambda_1  x. u $,  $(w \ v)
\ras u[x:=v]$.
\end{lemma}
\begin{proof}
By induction on $u$.
\end{proof}

\begin{lemma}\label{local}
The reduction $\tr$ is locally confluent on  fair terms.
\end{lemma}
\begin{proof}
 The only critical pairs are  the following.
 \begin{itemize}
   \item $t=(w \ u_3)$, $w \in \lambda_1  x.(u_1 \  u_2)$,   $t \tr t_1=(w_1 \ u_3 \
   (w_2 \ u_3))$ for $w_j \in \lambda_1  x. u_j$, and $t \tr t_2 =(w' \ u_3)$ for $w' \in \lambda_1  x.v$
   and $(u_1 \ u_2) \tr v$.
Both $t_1$ and $t_2$ reduces to $v[x:=u_3]$.
\item $t=(S_{i+2} \ r_1 \ r_2)\in \lambda_i x.(u_1 \ u_2)$, $x \in u_1$, $x \not\in u_2$
(for example), for some $u_1 \tr v_1$ such that  $x \not\in v_1$,
$t \tr t_1= (K_i \ (v_1 \ u_2))$  and $t \tr t_2= (S_{i+2} \ (K_i
\ v_1) \ (K_i \ u_2)))$. But $t_2 \tr t_1$.
 \end{itemize}
\end{proof}

\begin{theorem}\label{confluent}
The reduction $\tr$ is  confluent on  fair terms.
\end{theorem}
\begin{proof}
By Lemma \ref{new} and \ref{local}.
\end{proof}

\subsection{Proof of theorem \ref{thm_aux}}
In this section I will still denote by $\tr$ the reduction on {\em
combinators} given by rules (1, 2, 4) of definition \ref{def_red}.

\begin{definition}
\begin{itemize}
  \item Let $u$ be a combinator. A labelling of $u$ is a function that associates
  to each occurrence of $S$ (resp. $K,I$) in $u$ some $S_i$  (resp. some $K_i, I_i$).
  \item If $L$ is a labelling of $u$,
  I still denote by $L(u)$ the term obtained by replacing in $u$
  the symbols $S$ (resp. $K,I$)  by
 $L(S)$ (resp. $L(K)$, $L(I)$).

   \item Let $u$ be a term. I denote by $\theta(u)$ the combinator defined  by the following rules.
   $\theta(x)= x$,   $\theta(S_i)=S$, $\theta(K_i)=K$, $\theta(I_i)=I$ and
$\theta((u \ v)= (\theta(u) \ \theta(v))$
\item Let $u$ be a combinator and $L, L'$ be  labelling of $u$. I say that $L'$ is an extension
   of $L$ if the following holds.

\begin{enumerate}
  \item For each $S$ in $u$,

   - either
   $L(S)=L'(S)$

   - or $L(S)=S_0$ and  $L'(S)=S_i$ for $i=1, 2$ or $ 3$

   - or $L(S)=S_2$ or $L(S)=S_1$ and $L'(S)=S_3$.
  \item For each $K$ in $u$, $L(K)=L'(K)$ or $L(K)=K_0$ and $L'(K)=K_1$.
   \item For each $I$ in $u$, $L(I)=L'(I)$ or $L(I)=I_0$ and
   $L'(I)=I_1$.
\end{enumerate}

\end{itemize}
\end{definition}

 A labelling of $u$ is thus a way of marking redexes in $u$. The function $\theta$ consists in un-marking  terms to
get combinators. Extending a labelling means allowing more redexes
to be reduced.

\begin{lemma}
Let $u$ be a combinator and $L$ be a labelling of $u$. If $L(u)
\tr v$ then $u \tr \theta(v)$.
\end{lemma}
\begin{proof}
Immediate.
\end{proof}

\begin{lemma}\label{carac4}

 Assume $t=L(\lambda x. r) \in F $ for some $L,r$. Then, there is an extension $L'$ of $L$ such that $L'(\lambda x. r)\in \lambda_0
  x.v$ for some $v \in F$.

\end{lemma}
\begin{proof}
First note that, for combinators, $\lambda x. r$ represents a
single term and thus having written $t=L(\lambda x. r)$ is not a
typo !

 $L'$ is obtained by iterating the following algorithm.

-  If $x$ does not occur in $r$, choose $L'=L$. Since $t=(L(K) \
L(r))$, by Lemma \ref{carac}, $L(K)$ must be $K_0$ and thus, by
Lemma \ref{carac1}, $L(r) \in F$.

 - If $r=x$, choose $L'=L$. The argument is similar.

-  If $r=(r_1 \ r_2)$. Then $\lambda x. r=(S \ \lambda x. r_1 \
\lambda x. r_2)$. By Lemma \ref{carac}, $L(S)$ must be either
$S_0$ or $S_2$.
\\ If $L(S)=S_2$, by Lemma \ref{carac1},
$t \in \phi(\lambda_0 x.v)$ for some $v \in F$ (the term $\phi(u)$
is defined in Definition \ref{def-phi}). Thus $L$ satisfies the
desired property since, by Lemma \ref{carac5}, $t$ must be in
$\lambda_0 x. \phi (v)$.
\\ If
 $L(S)=S_0$, then, by Lemma \ref{carac1},
$L(\lambda x. r_i) \in F $. Choose $L'(S)=S_2$ and  iterate  the
algorithm with $L(\lambda x. r_j)$ for $j=1, 2$.

\end{proof}

\begin{lemma}\label{extend}
Let $t$ be a combinator and $L$ be a labelling of $t$ such that
$L(t)$ is fair. Assume that $t \tr t'$. Then, there is an
extension $L'$ of $L$ such that $L'(t)$ is fair and $L'(t) \tr v$
for some $v$ such that $\theta(v)=t'$.
\end{lemma}
\begin{proof}
By induction on $ nb(L(t))$. Look at the last rule that has been
used to show that $L(t)$ is fair.

{\em Rule (3)} :   a redex in $w \in \lambda_0 x. u$ is either a
redex in $u$ (and the result follows immediately from the \ih) or
it is of the form $(S_2 \ (K_0 \ u_1) \ (K_0 \ u_2))
  \tr (K_0 \ (u_1 \ u_2 ))$ and thus already appear in $L(t)$.

{\em Rule (5)} :   a redex in $\phi_x(u,f)$ is either a redex in
$u$ or in some $f(a)$ or a redex already in $L(t)$ and the result
follows immediately from the \ih.

{\em Rule (2)} : then $t=(t_1 \ t_2)$ and $L(t_1), L(t_2)$ are
fair. If the reduced redex is either in $t_1$ or $t_2$, the result
follows immediately from the \ih. Otherwise it has been  created
by the application of $t_1$ to $t_2$. I will only look at the
cases where the reduced redex starts with some $S$. The case of
$K$ and $I$ are similar and much simpler. For sake of simplicity I
will define $L'$ by only mentioning the labels that are changed.
We distinguish the different possible redexes.

\medskip
\noindent (a)  $t_1= (S \ u \ v)$ and  $t'= (u \ t_2 \ (v \ t_2))
$.

-  If $L(S)=S_0$ then, setting $L'(S)=S_1$ gives the desired
properties since, by Lemma \ref{carac1},
  $L(u), L(v)$ are in $F$ and thus $L'(t)$ also is in $F$.

-  $L(S)$ may not be $S_1$ or $S_3$ since, by Lemma \ref{carac},
it would have at least 3 arguments.

- If $L(S)=S_2$ then, by Lemma \ref{carac1}, $L(t)=\phi(w)$ for
some $w \in \lambda_0 x. v$ and some $v \in
   F$. Then, choosing $L'$ in such a way that $L'(t)=\phi(w_1 )$ for $w_1 \in \lambda_1  x.  v$ will give the desired properties .

\medskip
\noindent (b)
 $t_1= (S \ (K \ u)$, $t_2= (K \ v)$ and
$t'=(K \ (u \ v))$. Then $L(S)$ must be $S_0$ because otherwise,
by Lemma \ref{carac}, $S$ would have at least two arguments.
Similarly, we must have $L(K)=K_0$. Then, by Lemma \ref{carac1},
  $u,v$ are fair and thus setting
$L'(S)=S_2$ and $L'(K)=K_0$ gives the desired properties.

\ignore{\medskip \noindent (c) $t_1= (S \ (K \ u))$, $t_2=I$ and
$t'=u$. For the same reasons as in the previous case, we must have
$L(S)=S_0$, $L(K)=K_0$ and $L(I)=I_0$. Thus $u$ is fair and
setting $L'(S)=S_2$, $L'(K)=K_0$ and $L'(I)=I_0$ gives the desired
properties.}

\medskip
\noindent (c) $t_1=(S \ w_1)$ for $w_1 \in \lambda x.u_1$, $t_2\in
\lambda x.u_2$ and $t'\in \lambda x.v$ where $v$ is a reduct of
$(u_1 \ u_2)$. Again by Lemma \ref{carac},  we must have
$L(S)=S_0$. By Lemma \ref{carac1}, $L(w_1) \in F$. By Lemma
\ref{carac4}, extend $L$ so that $L'(u_i) \in F$.  Then setting
$L''$ in such a way that $L''(t)\in \lambda_0 x.(u_1 \ u_2)$ gives
the desired properties.

\medskip

{\em Rule (4)} : then $t=(S \ u_1 \ u_2 \ u_3)$, $L(S)=S_1$  and
the $L(u_i)$ are fair. If $t'= (u_1 \ u_3 \ (u_2 \ u_3))$ or if
the reduced redex is in some $u_i$ the result is trivial.
Otherwise this means that, for $i=1,2$ $u_i \in \lambda x. v_i$
and $t'= (w \ u_3)$ for some $w \in \lambda  x.v$ such that $v$ is
a reduct  of $(v_1 \ v_2)$. Then, by Lemma  \ref{carac4}, extend
$L$ so that $L'(v_i) \in F$ and choose $L''$ in such a way that
$L''(t)=(w'\ u_3)$ for $w' \in \lambda_1 x.(v_1 \ v_2) $.
\end{proof}

\begin{lemma}\label{conf}
Let $t$ be a combinator. Assume that  $t \tr v$ and $t \tr^* u$.
Then, there is a labelling $L$ of $u$ and a term $w$ such that
$L(u)$ is fair, $L(u) \ras w$ and $v \tr^* \theta(w)$.

\end{lemma}
\begin{proof}
By induction on the length $n$ of the reduction $t \tr^* u$.

\begin{itemize}
  \item If $n=1$, let $L_0$ be the labelling of $t$ obtained by indexing all
  the occurrences of $S,K$ and $I$ by 0. $L_0(t)$ is clearly fair. Apply
  Lemma \ref{extend} to  $t$, $L_0$ and the reduction $t \tr v$. This gives an extension $L_1$
  of $L_0$. Applying Lemma \ref{extend} to  $t$, $L_1$ and the reduction $t \tr u$ we get
  an extension $L_2$
  of $L_1$.
  Applying the confluence of $\tr$ on fair terms  (Theorem \ref{confluent}) to $L_2(t)$ gives the desired result.
  \item Otherwise, let $t \tr^* u_1 \tr u$. By the \ih,  let $L_1$ be a labelling of $u_1$ and
   $w_1$ be a term
such that $L_1(u_1)$ is fair, $L_1(u_1) \ras w_1$ and $v \tr
\theta(w_1)$. By Lemma \ref{extend}, let $L$ be a labelling of
$u_1$ that is an extension of $L$ such that $L(u_1)$ is fair and
$L(u_1)\tr r$ for $r$ such that $\theta(r)=u$. By theorem
\ref{confluent}, let $w$ be such that $r \ras w$ and $w_1 \ras w$.
Then $L,w$ have the desired properties.
 \end{itemize}
\end{proof}

\begin{proposition}\label{main1}
The reduction given by rules (1, 2, 4) of definition \ref{def_red}
is confluent.
\end{proposition}
\begin{proof}
It is enough to show that, if $t \tr u$ and $t \tr^* v$ then $u
\tr^* w$ and $v \tr^* w$ for some $w$. This follows immediately
from Lemma \ref{conf}.
\end{proof}

\begin{definition}
I denote by  $\supset$ the reduction defined by the following
rules.
\begin{enumerate}
  \item $(S \ (K \ u) \ I) \supset u$ \hspace{.5cm}  $(K \ u \ v)
  \supset u$ \hspace{.5cm} $(I \ u) \supset u$
  \item  $\lambda x.u
\supset \lambda x.v$ if $u \supset v$

\end{enumerate}
\end{definition}

\begin{lemma}\label{prepa3}
The reduction $\supset$ is confluent and commutes with $\tr$.

\end{lemma}
\begin{proof}
The reduction $\supset$ is strongly normalizing since it decreases
the size. Thus to prove the confluence, it is  thus enough to show
the local confluence and this is straightforward.
 Since  $\supset$ is also non duplicating, to prove the commutation with $\tr$, it is  enough to show the local
commutation and this is again straightforward. Note that the
reductions $(K \ u \ v)
  \supset u$, $(I \ u) \supset u$ that are already present in
  $\tr$ are used here to ensure the confluence of the only
  critical pair i.e. $(S \ (K \ u) \ I \ w) \supset (u \ w)$ and
  $(S \ (K \ u) \ I \ w)  \tr (K \ u \ w \ (I \ w))$.

\end{proof}

\noindent {\bf Theorem \ref{thm_aux}} {\em The reduction  given by
rules (1, 2, 3, 4) of definition \ref{def_red}  is confluent.}

\begin{proof}
Since $\rightarrow$ is the union of $\tr$ and $\supset$, the
result follows immediately from proposition \ref{main1} and Lemma
\ref{prepa3}.
\end{proof}

\section{Proof of theorem \ref{thm}}\label{s4}

\begin{definition}
I denote by $\vdash$ the reduction defined by the following rules.
\begin{enumerate}
  \item $(S \ (K \ u) \ I) \vdash u$ \hspace{.5cm}  $(K \ u \ v)
  \vdash u$ \hspace{.5cm} $(I \ u) \vdash u$
  \item  $[x]u
\vdash [x]v$ if $u \vdash v$

\end{enumerate}
\end{definition}

\begin{lemma}\label{prepa}
 The reduction $\vdash$ is confluent and commutes with $\succ$.
\end{lemma}
\begin{proof}
As in Lemma \ref{prepa3}
\end{proof}

\begin{lemma}\label{prepa2}
If $u \rightarrow^* v$ then $u \succ^* w$, $v \vdash^* w$ for some
$w$.
\end{lemma}
\begin{proof}
By induction on the length of the reduction $u \rightarrow^* v$.
Assume $u \rightarrow u_1 \rightarrow^* v$. If the level of the
reduction $u \rightarrow u_1$ is 0, the result follows immediately
from the \ih\ since then we also have $u \succ u_1$. Otherwise,
the reduction looks like $u=C[\lambda x. t] \rightarrow u_1=
C[\lambda x. t_1] \rightarrow^* v$ where $t \rightarrow t_1$. By
the \ih, we have $t \succ^* w_1$, $t_1 \vdash^* w_1$ for some
$w_1$ and thus $u \succ^* w_2$, $u_1 \vdash^* w_2$ where
$w_2=C[w_1]$. By the \ih\ we also have $u_1 \succ^* w$, $v
\vdash^* w$ for some $w$. By Lemma \ref{prepa}, we have $w_2
\succ^* w_3$ and $w \vdash^* w_3$ for some $w_3$ which is the term
we are looking for.
\end{proof}

\noindent {\bf Theorem \ref{thm}} {\em The reduction $\succ$ is
confluent}.

\medskip
\begin{proof}
Assume $t \succ^*t_1$ and $t \succ^*t_2$. Then $t_1 \equiv t_2$
and thus, by theorem \ref{equiv}, $t_1 \approx t_2$. Since
$\rightarrow$ is confluent we thus have $t_1 \rightarrow^* t_3$,
$t_2 \rightarrow^* t_3$ for some $t_3$. By Lemma  \ref{prepa2},
let $v_i$ be such that $t_i \succ^* v_i$ and $t_3 \vdash^* v_i$.
Since $\vdash$ is confluent, let $t_3$ be such that $v_1 \vdash
t_3$ and $v_2 \vdash t_3$. Since  $\vdash$ is clearly a
restriction of $\succ$, we have $t_i \succ^*t_3$
\end{proof}

\section{A standardization theorem}\label{sst}
In this section I prove a standardization theorem for the system
of section \ref{s3}. I study this system instead of the one of
section \ref{s2} because, as already mentioned in section
\ref{s2.2},  in the original system,  what could be the leftmost
redex is not clear at all.

Note that the following definition of a standard reduction does
not need the definition of the residue of a redex. It is a
definition by induction on $\langle lg(t \f t' ), size(t)\rangle$
where $lg(t \f t' )$ is the number of steps of the reduction. It
uses the idea that is implicit in \cite{david} and simply says
that a standard  reduction either reduces the head redex at the
first step or is not allowed to reduce it.
\begin{definition}
A reduction $t \f^* t'$ is standard ($t \fst t'$ for short) if it
satisfies  the following properties.
\begin{enumerate}
  \item $t=(x \ \overrightarrow{u})$, $t'=(x \
  \overrightarrow{u'})$ and, for each $i$, $u_i \fst u'_i$
  \item $t=(K \ u)$, $t'=(K \ u')$  and $u \fst u'$.
\item  $t=(S \ u)$, $t'=(S\ u')$ and $u \fst u'$.

  \item $t=(S \ u \ v)$ and

\begin{itemize}
  \item either $t'=(S \ u' \ v')$ for $u \fst
  u'$  and $v \fst v'$
  \item or the reduction is $t \f t_1 \ ... \ \f t_k \fst t'$ for some $k \geq 0$ such that
  $t_i=(S \ u_i \
  v_i)$, $u \fst
  u_k$, $v \fst
  v_k$ and

  - either $t_k=[x]w$,  $t'=[x]w'$, $w \fst w'$ and, for each $i<k$,
  $t_i$ cannot be written as $[x]r$ for some $r$

  - or $u_k=(K \ u'_k), v_k=(K \ v'_k)$, the reduction $t_k \fst
  t'$ is $t_k \f (K \ (u'_k \ v'_k)) \fst t'$ and, for each $i<k$,
  $t_i$ cannot be written as $(S \ (K \ u'_i) \ (K \ v'_i))$

  - or $u_k=(K \ u'_k), v_k=I$, the reduction $t_k \fst
  t'$ is $t_k \f u'_k  \fst t'$ and, for each $i<k$,
  $t_i$ cannot be written as $(S \ (K \ u'_i) \ I)$

\end{itemize}

  \item $t=(I \ u_1 \ ... \ u_n)$ for $n \geq 1$ and

\begin{itemize}
  \item either $t'=(I \ u'_1 \ ... \
  u'_n)$ for $u_i \fst u'_i$
  \item or the reduction is $t \f (u_1 \ ... \ u_n) \fst t'$

\end{itemize}
  \item $t=(K \ u_1 \ ... \ u_n)$ for $n \geq 2$ and

\begin{itemize}
  \item  either $t'=(K \ u'_1 \ ... \
  u'_n)$ for $u_i \fst u'_i$
  \item or the reduction is $t \f (u_1 \ u_3 \ ... \ u_n) \fst t'$

\end{itemize}
\item $t=(S \ u_1 \ ... \ u_n)$ for $n \geq 3$ and

\begin{itemize}
  \item  either $t'=(r \  u'_3 \ ... \
  u'_n)$ where $(S \ u_1 \  u_2) \fst r$ and  $u_i \fst u'_i$ for
  $i \geq 3$
  \item or the reduction is $t \f (u_1 \ u_3 \ ( u_2 \ u_3) \ u_4 \ ... \ u_n) \fst t'$
\end{itemize}
\end{enumerate}
\end{definition}

\begin{lemma}\label{prepa_st}
\begin{itemize}
  \item Assume $u_i \fst u'_i$ for each $i$. Then $(u_1  \ ... \
  u_n) \fst (u'_1  \ ... \
  u'_n)$
  \item Assume $u \fst [x]u'$. Then $(u \ v) \fst u'[x:=v]$

\end{itemize}
\end{lemma}
\begin{proof}
Easy.
\end{proof}

\begin{theorem}\label{st}
If $t \f^* t'$ then $t\fst t'$.
\end{theorem}
\begin{proof}
By induction on $lg(t \f^* t' )$. It is enough to show that if $t
\fst t' \f t''$ then $t \fst t''$. This is done by induction on
$\langle lg(t \fst t' ), size(t)\rangle$ and by case analysis. We
look at the rule that has been used to show $t \fst t'$ and then
what is the reduced redex in $t' \f t''$. I just consider two
cases. The first one is typical and easy. The second one is
similar but a bit more complex.

\begin{itemize}
  \item $t=(K \ u_1 \
... \ u_n)$ for $n \geq 2$.

\begin{itemize}
  \item If the reduction is $t \f (u_1 \ u_3 \ ... \ u_n) \fst t'$ we
apply the \ih\ to $(u_1 \ u_3 \ ... \ u_n)$ $ \fst t' \f t''$ to
get $(u_1 \ u_3 \ ... \ u_n) \fst t''$ and thus $t \f (u_1 \ u_3 \
... \ u_n) \fst t''$ is standard.
  \item  If the reduction is such that $t'
=(K \ u'_1 \ ... \
  u'_n)$ for $u_i \fst u'_i$ then

-  either $t''=(K \ u'_1 \ ... \
  u''_i \ ... \
  u'_n)$ for $u'_i \f u''_i$ and we apply the \ihb to $u_i \fst
  u'_i \f u''_i$ to get the result

  - or $t''=(u'_1 \ u'_3 \ ... \
  u'_n)$ and then $t \f (u_1 \ u_3 \
... \ u_n) \f^* (u'_1 \ u'_3 \ ... \
  u'_n)$ is standard by Lemma \ref{prepa_st}.
\end{itemize}
\item  $t=(S \ u_1 \ ... \ u_n)$ for $n \geq 3$ and
 $t'=(r \  u'_3 \ ... \
  u'_n)$ where $(S \ u_1 \  u_2) \fst r$ and  $u_i \fst u'_i$ for
  $i \geq 3$. Assume also  that $r=[x]a$, $x \notin r$ and $t''= (a \  u'_4 \ ... \
  u'_n)$. This means that, for $i=1,2$,  $u_i \fst [x]v_i$ and
  that $(v_1 \ v_2) \fst a$. But then, by Lemma \ref{prepa_st},
  $(u_i \ u_3) \fst v_i[x:=u_3]$. Thus, the following reduction
  is standard. $t \f (u_1 \ u_3 \ ( u_2 \ u_3) \ u_4 \ ... \ u_n)
  \fst (v_1[x:=u_3] \ v_2[x:=u_3] \  u_4 \ ... \ u_n) \fst
  (a[x:=u_3] \  u_4 \ ... \ u_n) \fst (a \  u'_4 \ ... \ u'_n)=  t''$.

\end{itemize}
\end{proof}

\ignore{

The next result will not be used later but may be also of some
interest. Note that, even though its proof looks like the one of
the standardization theorem, it is not a corollary of it.

\begin{definition}
An infinite reduction $t \f t_1 \ t_2 \f ...$ is standard i, for
each $i$, the reduction $t \f^* t_i$ is standard.
\end{definition}

\begin{theorem}\label{stbis}
If $t$ is not strongly normalizable then $t$ has an infinite
standard reduction.
\end{theorem}
\begin{proof}
Assuming $t \notin SN$,  an infinite standard reduction is given
by applying the following algorithm.
\begin{enumerate}
  \item If $t=(x \ \overrightarrow{u})$,  then some $u_i$ is not
  in $SN$. Apply the algorithm to such an $u_i$
  \item If $t=(K \ u)$, then $u \notin SN$. Apply the algorithm to
  $u$

  \item If $t=(S \ u \ v)$. If $u \notin SN$ (resp. $v \notin
  SN$), apply the algorithm to $u$ (resp.$v$). Otherwise, this
  means that $u \f^* [x]u_1$ and $v \f^* [x]v_1$ for some $u_1,v_1$ such that $(u_1 \ v_1)
  \notin SN$. Then, by theorem \ref{st}, $u \fst [x]u_1$ and $v \fst
  [x]v_1$. Apply the algorithm to $(u_1 \ v_1)$.

  \item $t=(I \ u_1 \ ... \ u_n)$ for $n \geq 1$. If some $u_i$ is
  not in $SN$, apply the algorithm to such an $u_i$. Otherwise, apply the algorithm to
  $( u_1 \ ... \ u_n)$.

  \item $t=(K \ u_1 \ ... \ u_n)$ for $n \geq 2$. If some $u_i$ is
  not in $SN$, apply the algorithm to such an $u_i$. Otherwise, apply the algorithm to
  $( u_1 \ u_3 \ ... \ u_n)$.
\item $t=(S \ u_1 \ ... \ u_n)$ for $n \geq 3$.  If some $u_i$ is
  not in $SN$, apply the algorithm to such an $u_i$. If $t=(S \ u_1  u_2) \notin
  SN$ apply the algorithm to $(S \ u_1 \ u_2)$. Otherwise use a
  standard reduction to get $(S \ u_1  u_2) \fst r$, reduce the
  head redex of $(r \ u_3 \ ... \ u_n)$ to get $t_1$ and apply the
  algorithm to $t_1$.
\end{enumerate}
It follows easily from theorem \ref{st} that the given algorithm
produce an infinite standard reduction.
\end{proof}

}

\section{Strong normalization of the typed calculus}\label{ssn}

In this section I prove the strong normalization of the auxiliary
system of section \ref{s3}. Note that the system of section
\ref{s2} is not strongly normalizing even though this is for the
following bad reason. Let $t=(S \ x \ x)$. Then $t=[y](S \ x \ x \
y ) \succ [y](x \ y \ (x \ y))=t$.

The types are the simple types i.e. constructed from basic types
with the arrow. The typing rules are the usual ones i.e.   $I$ has
type $A \f A$, $K$ has type $A \f B \f C$, $S$ has type $(A \f B
\f C) \f (A \f B) \f A \f C$ for every types $A,B,C$ and, finally,
if $u$ has type $A \f B$ and $v$ has type $A$ then $(u \ v)$ has
type $B$.

\begin{definition}
\begin{itemize}
  \item A combinator $t$ is highly normalizing ($t \in HN$ for short) if
it can be obtained by the following rules.
\begin{enumerate}
  \item $t=S$ or $t=K$ or $t=I$ or $t=(x \ t_1 \ ... \ t_n)$ for
  $t_1,  \ ...,  \ t_n \in HN$.
  \item $t=(K \ t_1)$ or $t=(S \ t_1)$  for $t_1 \in
  HN$
  \item $t=(S \ t_1 \ t_2)$ for $(t_1 \ x \ (t_2 \ x)) \in HN$ where
  $x$ is a variable.
  \item $t=(I \  t_1 \ ... \ t_n)$ for $n \geq 1$ and $(t_1 \  t_2 \ ... \ t_n) \in HN$

\item $t=(K \  t_1 \ ... \ t_n)$ for $n \geq 2$, $(t_1 \  t_3 \ ... \ t_n) \in
HN$and $t_2 \in HN$
\item $t=(S \  t_1 \ ... \ t_n)$ for $n \geq 3$ and $(t_1 \ t_3 \ ( t_2 \ t_3) \ t_4 \ ... \ t_n) \in HN$
\end{enumerate}

  \item If $t \in HN$ we denote by $\eta(t)$ the number of rules
  that have been used to show $t \in HN$.
\end{itemize}
\end{definition}

We have introduced this  notion of normalization which is stronger
than the usual one (see the next Lemma) because the proof of Lemma
\ref{main_sn} below would not work if $HN$ was replaced by $SN$.

\begin{lemma}
If $t \in HN$  then $t$ is strongly normalizing.
\end{lemma}

\begin{proof}
By induction on $\eta(t)$. The non trivial cases are when the last
rule that has been applied to prove $t \in HN$ is (3) or (6).
\begin{itemize}
  \item Assume first $t=(S \ t_1 \ t_2)$. Then, by the \ih, $t'=(t_1 \ x \ (t_2 \ x)) \in
  SN$ and thus $t_1, t_2 \in SN$. Thus an infinite reduction of
  $t$ must look like
 $t \f^* t'' \f^*  ...$ where for some $v_i$,  $t_i \f^* \lambda x.v_i$ and

 - either the reduction of $t''$ is in  $(v_1 \ v_2)$.
  But $(v_1 \ v_2) \in SN$
  since  $t' \in SN$ and $t' \f^* (\lambda x. v_1 \ x \ (\lambda x.v_2 \ x))
  \f^* (v_1 \ v_2)$. Contradiction.

  - or $\lambda x.v_i=(K \ v_i)$ and the reduction  is
  $t''=(S \ (K \ v_1) \ (K \ v_2) \f (K (v_1 \ v_2)) \f^* ...$
  This is  impossible since $t' \in SN$ and $t' \f^* (v_1 \
  v_2)$.

  -or $\lambda x.v_1=(K \ v_1)$, $\lambda x.v_2=I$  and the reduction is
  $t''=(S \ (K \ v_1) \ I \f v_1 \f^* ...$ This is impossible since $t' \in SN$ and $t' \f^* v_1$.
  \item $t=(S \  t_1 \ ... \ t_n)$. Again, by the \ih, $t'= (t_1 \ t_3 \ ( t_2 \ t_3) \ t_4 \ ... \ t_n) \in
  SN$. Thus the $t_i$ are in $SN$ and also $(S \  t_1 \  t_2) \in SN$. The first point is clear. For the second,
   we argue as follows. Reasoning as in the previous  case, it is
   enough to show that $(t_1 \ x \ (t_2 \ x)) \in
  SN$. If it was not the case then $(t_1 \ t_3 \ ( t_2 \ t_3)$
  would also not been in $SN$ and this contradicts the fact that
  $t' \in SN$.
  Thus an infinite reduction of $t$
  looks like
 $t \f^* (r \  t'_3 \ ... \ t'_n)  \f t'' \f^*
  ...$ where $r$ is a reduct of $(S \ t_1 \ t_2)$ and $t''$ is obtained by an interaction between $r$ and its arguments.
  But we have shown (in the proof of theorem \ref{st}) that then
  $t'$ reduces to $t"$ and this is a contradiction.
\end{itemize}
\end{proof}

\begin{lemma}\label{main_sn}
Let $t$ be a combinator and $\sig$ be a substitution such that all
the variables in the domain of $\sig$ have the same type. Assume
$t \in HN$ and the image of $\sig$ is included in $HN$. Then
$\sig(t) \in HN$.
\end{lemma}

\begin{proof}
By induction on $\langle type(\sig), \eta(t)\rangle$. Look at the
last rule that has been used to prove $t \in HN$. The only non
trivial case is when $t=(x \ t_1 \ ... \ t_n)$ and $x \in
dom(\sig)$. By the \ih,  $u_i=\sig(t_i) \in HN$. We now have to
distinguish the different possible values for $\sig(x)$. The most
difficult case (the other ones are similar or trivial) is when
$\sig(x)=(S \ a_1 \ a_2)$. We have to show that $t'=(a_1 \ u_1 \ (
a_2 \ u_1) \ u_2 \ ... \ u_n) \in
  HN$. But $t'=\tau((z \  u_2 \ ... \ u_n))$ where $z$ is a fresh
  variable such that $\tau(z)=(a_1 \
u_1 \ ( a_2 \ u_1))$. But $type(z) < type(x)$ and, by the \ih, it
is thus enough to show that $t''=(a_1 \ u_1 \ ( a_2 \ u_1)) \in
HN$. But $t''=\tau'((a_1 \ z' \ ( a_2 \ z')))$ where $z'$ is a
fresh
  variable such that $\tau'(z')=u_1$. Since $type(z') < type(x)$
  and $(a_1 \ z' \ ( a_2 \ z')) \in HN$ (because $(S \ a_1 \ a_2)\in HN$),
  the result follows from the \ih.
\end{proof}

\begin{corollary}
Every typed combinator $t$ is in $HN$ and thus in $SN$.
\end{corollary}
\begin{proof}
By induction on the size of $t$ using $(u \ v)=(x \ v)[x:=u]$ and
Lemma \ref{main_sn}.
\end{proof}

\section{Final remarks}\label{s5}

Though intuitively quite simple, the given proof of confluence is
technically rather involved and, in particular, it is more
elaborate than the one using the confluence of the $\l$-calculus.
Thus, one may wonder about the real use of such a proof even if
this is the condition to have a self contained theory.
 I will
argue for another reason.

 Combinatory Logic
 somehow looks like a calculus with explicit substitutions. Though $([x]u \ v)$ is not exactly
 the {\em explicit} substitution $u[x:=v]$, it
 has often to be understood in this way. In particular, the reduction
  $([x](u_1 \ u_2) \ v) \rightarrow ([x] u_1 \ v \ ([x]u_2 \ v))$
 looks like the propagation of the substitution into the two branches of the application.
But proving confluence for such calculi
 is usually not trivial simply because the usual methods (parallel reductions or finite developments) need
 definitions that are not clear.

I thus hope that the given proof will help in finding simple
proofs for calculi with explicit substitutions.

\medskip

\noindent {\bf Acknowledgments}

I wish to thank R. Hindley and the anonymous referee for helpful
comments on previous versions of this paper.

\medskip

\noindent {\bf Added in proof}

Shortly after sending the first version of this paper, I have been
informed by R.Hindley and P.Minari that the later has also written
(more or less at the same time)  a direct proof of the confluence
of combinatory strong reduction. This proof is completely
different from the one given here. See the TLCA list of open
problem or \cite{minari}.

\end{document}